\documentclass[10pt]{amsart}
\usepackage{bbm}
\usepackage[all]{xy}

\theoremstyle{definition}

\theoremstyle{remark}

\numberwithin{equation}{section}

\def\Ext{\mathop\mathrm{Ext}\nolimits}

\def\la{\lambda}

\def\aK{\mathbbm{k}}	
\def\fW{\mathbf{w}}		\def\fE{\mathbf{e}}

\def\qR{{\boldsymbol R}}		\def\qS{{\boldsymbol S}}

\def\kE{\mathcal{E}} 	\def\kF{\mathcal{F}}

\def\+{\oplus}
\def\*{\otimes}
\def\8{\infty}

\def\ra{\rightarrow}

\def\xx{\times}

\def\mtr#1{\begin{pmatrix}#1\end{pmatrix}}

 \newcommand{\hhh}{\rule{0pt}{9pt}}

\def\cm{Cohen--Macaulay}

\def\ti{\tilde}

\def\b{\bar}

\begin{document}

\title[Cohen-Macaulay modules over $T_{36}$]{Cohen-Macaulay modules over the plane curve singularity of type $T_{36}$}

\author[Y. Drozd]{Yuriy A. Drozd}
\address{Institute of Mathematics, National Academy of Sciences, 01601 Kyiv, Ukraine}
\email{y.a.drozd@gmail.com,\,drozd@imath.kiev.ua}
\urladdr{http://www.imath.kiev.ua/$\sim$drozd}

\author[O. Tovpyha]{Oleksii Tovpyha}
\address{Institute of Mathematics, National Academy of Sciences, 01601 Kyiv, Ukraine}
\email{tovpyha@gmail.com}

\subjclass[2010]{13C14, 14H20}
\keywords{Cohen--Macaulay modules, matrix factorizations, bimodule problems, bunches of chains}

\begin{abstract}
For a wide class of Cohen--Macaulay modules over the local ring of the plane curve singularity of type $T_{36}$ we describe explicitly the 
corresponding matrix factorizations. The calculations are based on the technique of matrix problems, in particular, representations of bunches 
of chains.
\end{abstract}
\maketitle

\section{Introduction}
\label{s1} 

Let $\aK$ be an algebraically closed field, $\qS=\aK[[x,y]]$. Recall that the complete local ring of the plane curve singularity of type $T_{36}$ is $\qR=\qS/(F)$, where $F=x(x-y^2)(x-\la y^2)$ and $\la\in\aK\setminus\{0,1\}$. In this paper we present explicit description of a wide class of maximal
 \cm\ modules over the ring $\qR$ called \emph{modules of the first level}. Note that $T_{36}$ is one of the critical singularities of tame 
 Cohen-Macaulay representation type \cite{dg}. Till now, only for the singularities of type $T_{44}$ matrix factorizations have been described
 \cite{dt1,dt2}.
 
\section{Matrix problem and the first reduction}
\label{s2} 

 So, let $\qR=\aK[[x,y]]/(F)$, where $F=x(x-y^2)(x-\la y^2)\ (\la\in\aK\setminus\{0,1\})$.  
We consider $\qR$ as the subring of the direct product $\ti{\qR} = \qR_1\xx\qR_2\xx\qR_3$, where all $\qR_i=\aK[[t]]$, generated by the elements $x=(0, t^2, \la t^2)$ and $y=(t,t,t)$. We denote by $\qR_{12}$ the projection of $\qR$ onto $\qR_1\xx\qR_2$. It is generated by $(t,t,0)$ and $(0,t^2,0)$ and $\qR_{12} \simeq \aK[[x,y]]/x(x-y^2)$. It is a singularity of type $A_2$, so all indecomposable $\qR_{12}$-modules are $\qR_1,\,\qR_2,\, \qR_{12}$ and $\qR'_{12} = \qR_{12}[t_1]$, where $t_1 =(t,0,0)$. We also set $t_2=(0,t,0)$.

Denote $\ti{M} = \ti{\qR}M = M_1 \+ M_2 \+ M_3$, where $M_i$ is an $\qR_i$-module. There is an exact sequence:
 $$0 \ra M_{12} \ra M \ra M_3 \ra 0,$$ where $M_{12}=M\cap(M_1\+M_2)$ is an $\qR_{12}$-module, hence 
 $M_{12} = m_1\qR_1 \+ m_2\qR_2 \+ m_{12}\qR_{12} \+ m_{12}'\qR_{12}'$. So $M$ gives an element $\xi \in \Ext^1_\qR(M_{3}, M_{12})$. 
 There is an exact sequence
 $$0 \ra (x - \la y^2) \qR \ra \qR \ra \qR_{3} \ra 0,$$
whence
 $$\Ext^1_\qR(M_{3}, M_{12}) = M_{12}/(x - \la y^2)M_{12}$$

In the table nearby we present bases of the modules $\Ext^1_\qR(\qR_3, N)$, where $N \in {\qR_1, \qR_2, \qR_{12}, \qR'_{12}}$. In this table each
element actually denotes its residue class modulo $x-\la y^2$, $1_s\ (s\in\{1,2,12\})$ is the identity element of $\qR_s$, $t_{12}=t_1+t_2$ 
and $\mu = (1-\la)/\la$.
 \[
  \begin{array}{|c|c|}
  \hline 
  \hhh \qR_1 & 1_1,\ t 1_1 \\
  \hline
  \hhh \qR_2 & 1_2,\ t_2 \\
  \hline
  \hhh \qR'_{12} & 1_{12},\ t_1,\ t_2,\ t^2_1=\mu t^2_2 \\
  \hline
  \hhh \qR_{12}  & 1_{12},\ t_{12},\ t^2_1=\mu t^2_2,\ t^3_1=\mu t^3_2\\
  \hline
  \end{array}
 \] 
 
Homomorphisms $\qR_{12} \ra \qR_i$ induce the maps $\Ext^1_\qR(\qR_3, \qR_{12}) \ra \Ext^1_\qR(\qR_3, \qR_i)$ which map $1_{12} \mapsto 1_i$ and $t_{12} \mapsto t_i$. The embedding $\qR'_{12} \ra \qR_{12}$ which maps $1_{12} \mapsto t_{12}$ induces the map 
$\Ext^1_\qR(\qR_3, \qR'_{12}) \ra \Ext^1_\qR(\qR_3, \qR_{12})$ that coincide with the multiplication by $t_{12}$. The embeddings $\qR_i \ra \qR'_{12}$ such that $1_i \mapsto t_i$ induce the maps $\Ext^1_\qR(\qR_3, \qR_i) \ra \Ext^1_\qR(\qR_3, \qR'_{12})$ that coincide with the multiplication by $t_i$.

In particular, if we only consider free terms, we obtain representations of the partially ordered set of width 2:
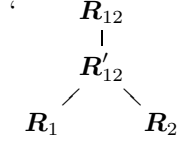
\[`
   \xymatrix@C=.3ex@R=1.3ex{ & \qR_{12} \ar@{-}[d] \\  & \qR'_{12} \ar@{-}[dl] \ar@{-}[dr] \\ \qR_1 && \qR_2 }
\]
(see \cite{nr}), hence the matrix $A_0$ of free terms can be reduced to the form:
$$
A_0=\left(
\begin{array}{c|c|c|c|c|c} 
	0 & 0   & 0   & I_1 & 0 & 0 \\
	0 & I_1 & 0   & 0   & 0 & 0 \\
	0 & 0   & 0   & 0   & 0 & 0 \\\hline
	0 & 0   & I_2 & 0   & 0 & 0 \\
	0 & I_2 & 0   & 0   & 0 & 0 \\ 
	0 & 0   & 0   & 0   & 0 & 0 \\\hline
	0 & 0   & 0   & 0   & I_{12} & 0 \\
	0 & 0   & 0   & 0   & 0 & 0 \\\hline
	0 & 0   & 0   & 0   & 0 & I_{12} \\
	0 & 0   & 0   & 0   & 0 & 0 \\	
\end{array}
\right),
$$
 where $I_s$ means $1_sI$ for the identity matrices $I$ of the appropriate size (maybe different for different matrices). 
 In what follows we always suppose that $A_0$ is of this form.

\section{Modules of the first level}
\label{s3} 

Let now $A = A_0+tA_1$ where $A_1$ is also divided into blocks in the same way as $A_0$. Using automorphisms of $M_3$ we can make zero the 1st, 4th, 7th and 9th rows of $A_1$, as well as one of the 2rd or 5th rows.

Using automorphisms of $M_{12}$ we can also make zero all columns in  $A_1$, except the 1st one and the parts of the 2nd, 3rd and 4th columns in the 10th row, where we can only delete all terms containing $t^2$ or $t^3$. In particular, the terms $1_{12}$ from the last vertical stripe become
direct summands of the whole matrix $A$. So in what follows we can omit this column.  We always suppose that $A_1$ has this form.

 Let $A_1^0$ be the free term of the matrix $A_1$. The non-zero part of its 10th row can be considered as  representation of the partially
 ordered set:
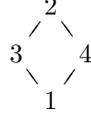
\[
   \xymatrix@C=.5ex@R=1.3ex{ & 2 \ar@{-}[dl] \ar@{-}[dr]  \\  3 \ar@{-}[dr] & & 4 \ar@{-}[dl] \\ & 1 } 
\]
 Hence it can be reduced to the form
 \[
 \left( \begin{array}{cc|cc|cc|cc}
  0&I &0&0&0&0&0&0 \\ 0&0&I&0&0&0&0&0 \\
 0&0&0&0&0& I& 0&I\\ 0&0&0&0&0&0&0&0
 \end{array}\right),
 \]
 where $I$ is again an identity matrix of the appropriate size (maybe different for different matrices). 
 Then the whole matrix $A$ modulo $t^2$ can be reduced to the form:
\[ 
 \begin{array}{c|cc|cc|cc|cc|c}
 &1&2& & &3& &4&5& \\\hline 
 &0&0&0&0&0&0&I_1&0&0\\  &0&0&0&0&0&0&0&I_1&0 \\
 1&A_{11}t^* & 0& I_1&0&0&0&0&0&0 \\  2& A_{21}t^*&0 & 0& I_1&0&0&0&0&0 \\ 
 3& A_{31}t_1&0&0&0&0&0&0&0 &0 \\\hline
 &0&0&0&0&I_2&0&0&0&0\\  &0&0&0&0&0&I_2&0&0&0 \\
  & 0&0 &  I_2&0&0&0&0&0&0 \\  & 0&0&0&I_2&0&0&0&0 &0 \\
 4& A_{41}t_2&0&0&0&0&0&0&0 &0 \\ \hline
 5&A_{51}t^* &A_{52}t^*  & 0&0&0&0&0&0& I_{12} \\
 6& A_{61} & A_{62}t^* & 0&0& A_{63}t_1 & 0 & A_{64}t_2 &A_{65}t^* &0 \\
 \hline
 &0&t_{12}I &0&0&0&0&0&0&0 \\ &0&0&t_{12}I&0&0&0&0&0&0 \\
 &0&0&0&0&0& t_{12}I& 0&t_{12}I& 0\\ &0&0&0&0&0&0&0&0&0 
\end{array} 
 \]
 Here the symbol $t^*$ means that in this block $t_1=-t_2$, and $A_{61}$ is a matrix pencil $X_1t_1+X_2t_2$. The horizontal lines show the
 division of $A$ into the stripes such that the 1st stripe corresponds to $R_1$, the 2nd to $R_2$, the 3rd to $R'_{12}$ and the 4th to $R_{12}$. 
 Moreover, as in the matrix $A$ we have $t_1^2=\mu t_2^2$ with $\mu\ne-1$, one can delete all terms with $t_i^2$ and $t_i^3$ everywhere 
 except the last block of the first column.
 
 The endomorphisms of $M_3$ and $M_{12}$ which do not destroy the shape of the matrices $A_0$ and $A_1^0$ induce the transfromations of 
 columns that can be described by the scheme
  
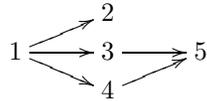
\[
   \xymatrix@R=.5ex{ & 2 \\ 1 \ar[ur] \ar[r]  \ar[dr] & 3 \ar[r] & 5 \\  & 4 \ar[ur]}
  \] 
  and the transformations of rows that can be described by the scheme
  
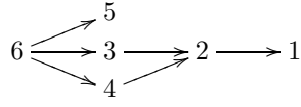
\[
   \xymatrix@R=.5ex{ & 5 \\ 6 \ar[ur] \ar[r]  \ar[dr] & 3 \ar[r] & 2 \ar[r] &1 \\  & 4 \ar[ur]}
  \]
  For the matrix $A_{61}$ it means that we can add the rows of $X_1$ to those of $A_{i1}$ for $i\in\{1,2,3\}$ and the rows of $X_2$ to the
  rows of $A_{i1}$ for $i\in\{1,2,4\}$. In  the same way, the columns of $X_1$ can be added to those of $A_{6j}$ for $j\in\{3,5\}$, while the columns
  of $X_2$ can be added to those of $A_{6j}$ for $j\in\{4,5\}$.
  
  The indecomposable matrix pencils (representations of the Kronecker quiver) are described in \cite{gan,ri}. In \cite{ri} the morphisms between 
  indecomposable representations are also described. It implies that the matrix $A_{61}$ is a direct sum of  the following matrices:
\begin{align*}
  & A(n)=\mtr{t_1&t_2&0&\dots&0\\0&t_1&t_2&\dots&0\\\hdotsfor5\\0&0&0&\dots&t_2\\0&0&0&\dots&t_1},\\
  & B(n)=\mtr{t_2&t_1&0&\dots&0\\0&t_2&t_1&\dots&0\\\hdotsfor5\\0&0&0&\dots&t_1\\0&0&0&\dots&t_2},\\
  & C(n)=\mtr{t_1&t_2&0&\dots&0&0\\ 0&t_1&t_2&\dots&0&0\\\hdotsfor6\\0&0&0&\dots&t_1&t_2}\\
  & D(n)=C(n)^\top\\
\end{align*}

It's easy to see that:
\begin{itemize}
	\item If $A_{61} = A(n)$ than we can make zero all matrices above $A_{61}$ except the 1st column of $A_{41}$, and all matrices to the right of 
	$A_{61}$ except the last row of $A_{64}$.
	\item If $A_{61} = B(n)$ than we can make zero all matrices above $A_{61}$ except the 1st column of $A_{31}$, and all matrices to the right 
	of $A_{61}$ except the last row of $A_{63}$.
	\item If $A_{61} = C(n)$ than we can make zero all matrices to the right of $A_{61}$, and all matrices above $A_{61}$ except the last column 
	of $A_{31}$, the 1st column of $A_{41}$ and one (any chosen) of the columns of the matrix $A_{51}$.
	\item If $A_{61} = D(n)$ than we can make zero all matrices above $A_{61}$, and all matrices to the right of $A_{61}$ except the last row of 
	$A_{63}$, the 1st row of $A_{64}$ and one (any chosen) of the rows of the matrix $A_{62}$.
\end{itemize}
Hence in the non-zero part of $A_{51}$ (and, respectively, $A_{62}$) we can left one non-zero element above each block of $C(n)$ (and, respectively, 
$D(n)$). Therefore, except the summands $A(n), B(n), C(n), D(n)$ in the blocks $A_{ij}$, ($i=5,6$, $j=1,2$), we will also have the summands of the
form $C'(n)$, with one additional element in $A_{51}$-part as compared to $C(n)$, and $D'(n)$, with one additional element in $A_{62}$-part as 
compared to $D(n)$). So we can suppose that $C'(n)$ looks like $B(n)^\top$, but with the 1st row from $A_{51}$, and $D'(n)$ looks like $A(n)^\top$, but 
with the last column from $A_{62}$. 

One can see now that now we can make zero all elements of the matrix $A_{52}$ except those which are in the zero rows of $A_{51}$ and zero ]columns of $A_{62}$. The remaining part of $A_{51}$ can be reduced to the form $$\mtr{I&0\\0&0}.$$ It gives  direct summands of the whole
matrix $A$ of the form $$\mtr{t_1\\\hline t_{12}}$$ (certainly, $t_1$ can be replace here by $t_2$). Therefore, in what follows we can suppose that 
$A_{52}=0$. Analogously, we can suppose that the matrices $A_{11},A_{12}$ and $A_{65}$ are also zero. Otherwise we obtain direct summands of
$A$. For instance, if $A_{65}\ne0$, all non-zero elements are in the rows which do not belong to the non-zero parts of $A_{66}$ and $A_{65}$.
So they give direct summands of the form
\[
 \mtr{0&1_1\\\hline 1_2&0\\\hline 0&t_2\\\hline t_{12}&t_{12}}.
\]

The description of homomorphisms between the representations of the Kronecker quuiver \cite{ri} show that we can add the non-zero columns
over $A(n)$ (respectively, $B(n)$\,) to those over $A(m)$ (respectively, $B(m)$\,) for $m>n$, and the same for the non-zero  rows to the right of $A(n)$ 
or $B(n)$. We can also add columns the non-zero columns over $C(n)$ (respectively, non-zero rows to the right of $D(n)$\,) to those of $C(m)$
(respectively, of $D(m)$), where $n<m$, as well as to those of $A(k)$ and $B(k)$ for any $k$. It means that the possible transformations of these
columns and rows can be considered as \emph{representations of a bunch of chains} in the sense of \cite{bo} or \cite[Appendix~B]{db} (we use the formulation of the second paper). Namely, we have the next pairs of chains:
\begin{itemize}
 \item  $\kE_1=\{a_i, d_i, d_i' \,|\, i \in \mathbb{N} \}$, $\kF_1=\{c_3\}$
 \item  $\kE_2=\{b_i, \tilde{d}_i, \tilde{d}_i' \,|\, i \in \mathbb{N} \}$, $\kF_2=\{c_4\}$
 \item  $\kE_3=\{r_3\}$, $\kF_3=\{\tilde{a}_i, c_i, c_i' \,|\, i \in \mathbb{N} \}$
 \item  $\kE_4=\{r_4\}$, $\kF_4=\{\tilde{b}_i,  \tilde{c}_i, \tilde{c}_i' \,|\, i \in \mathbb{N} \}$
\end{itemize}  
with the relation $\sim$: $$a_i \sim \tilde{a}_i,\ b_i \sim \tilde{b}_i,\, c_i \sim \tilde{c}_i,\ d_i \sim \tilde{d}_i,\ c'_i \sim \tilde{c}'_i,\ d'_i \sim \tilde{d}'_i\ 
(i \in \mathbb{N}).$$
Here $r_3, r_4$ corresponds to $A_{31}, A_{41}$ respectively and $c_3, c_4$ corresponds to $A_{63}, A_{64}$ respectively.

Now we use the description of the indecomposable representations of this bunch of chains from \cite{bo,db}. In our case they correspond to the following words in the alphabeth $\{a_i,\tilde{a}_i,b_i,\tilde{b}_i,c_i,\tilde{c}_i,d_i,\tilde{d}_i,c_i',\tilde{c}_i',d_i',\tilde{d}_i',c_3,r_4,c_4,r_3,-,\sim\}$:
\begin{itemize}
 \item 4 type of words with $a_i$, $i \in \mathbb{N}$: 
 
 $\fW_a(i) = r_3 - \ti{a}_i \sim a_i - c_3$ and 3 shorter words: $r_3 - \ti{a}_i \sim a_i$, $\ti{a}_i \sim a_i - c_3$, $\ti{a}_i \sim a_i$;
 
 \item 4 type of words with $b_i$, $i \in \mathbb{N}$:
 
 $\fW_b(i) = r_4 - \ti{b}_i \sim b_i - c_4$ and 3 shorter words: $r_4 - \ti{b}_i \sim b_i$, $\ti{b}_i \sim b_i - c_4$, $\ti{b}_i \sim b_i$;
 
 \item 4 type of words with $c_i$, $i \in \mathbb{N}$: 
 
 $\fW_c(i) = r_4 - \ti{c}_i \sim c_i - r_3$ and 3 shorter words: $r_4 - \ti{c}_i \sim c_i$, $\ti{c}_i \sim c_i - r_3$, $\ti{c}_i \sim c_i$;
 
 \item 4 type of words with $d_i$, $i \in \mathbb{N}$: 
 
 $\fW_d(i) = c_4 - \ti{d}_i \sim d_i - c_3$ and 3 shorter words: $c_4 - \ti{d}_i \sim d_i$, $\ti{d}_i \sim d_i - c_3$, $\ti{d}_i \sim d_i$;
 
 \item 4 type of words with $c_i'$, $i \in \mathbb{N}$:
 
 $\fW_c'(i) = r_4 - \ti{c}_i' \sim c_i' - r_3$ and 3 shorter words: $r_4 - \ti{c}_i' \sim c_i'$, $\ti{c}_i' \sim c_i' - r_3$, $\ti{c}_i' \sim c_i'$;
 
 \item 4 type of words with $d_i'$, $i \in \mathbb{N}$:
 
 $\fW_d'(i) = c_4 - \ti{d}_i' \sim d_i' - c_3$ and 3 shorter words: $c_4 - \ti{d}_i' \sim d_i'$, $\ti{d}_i' \sim d_i' - c_3$, $\ti{d}_i' \sim d_i'$;
\end{itemize} 

Following the construction of indecomposable representations from \cite{bo}, we construct the matrices corresponding to these words, see Table~1 below.  

$$P_a(n)=\left(\begin{array}{c|c} 0 & 1 \\ \hline t_2 \fE_1 & 0 \\  \hline \hhh A(n) & t_2 \fE_n^\top \end{array}\right)$$

$$P_b(n)=\left(\begin{array}{c|c} t_1 \fE_1 & 0 \\ \hline 0 & 1 \\ \hline \hhh B(n) & t_1 \fE_n^\top \end{array}\right)$$

$$P_c(n)=\left(\begin{array}{c} t_1 \fE_1 \\ \hline t_2 \fE_1 \\ \hline \hhh C(n) \end{array}\right)$$

$$P_d(n)=\left(\begin{array}{c|c|c} 0 & 1 & 0 \\ \hline 0 & 0 & 1 \\ \hline \hhh D(n) & t_2 \fE_{n+1}^\top & t_1 \fE_{n+1}^\top\end{array}\right)$$

$$P_c'(n)=\left(\begin{array}{c|c} t_1 \fE_1 & 0 \\ \hline t_2 \fE_1 & 0 \\ \hline \hhh C'(n) & t_{12} \fE_1^\top \end{array}\right)$$

$$P_d'(n)=\left(\begin{array}{c|c|c} 0 & 1 & 0 \\ \hline 0 & 0 & 1 \\ \hline \hhh D'(n) & t_2 \fE_{n+1}^\top & t_1 \fE_{n+1}^\top \\ \hline t_{12} \fE_{n+1} & 0 & 0 \end{array}\right)$$

Here $t_r = y 1_r$ and $\fE_n=(0,0,\dots,0,1)$, $\fE_1=(1,0,\dots,0)$ and $^\top$ means the transposition.

\section{Generators and relations. Example}
\label{s4} 
   
Now we calculate matrix factorizations of the polynomial $F=x(x-y^2)(x-\la y^2)$ corresponding to the indecomposable Cohen-Macaulay modules over $\qR$. In other words, we find minimal sets of generators for these modules and minimal sets of relations for these generators.

In order to make smaller the arising matrices, we denote $z=x-y^2$ an $z'=x-\la y^2$.  Thus $F=xzz'$.

We do detailed calculations for the word $\fW_a(2) = r_3 - \ti{a}_2 \sim a_2 - c_3$. Since all calculations are similar, for other words we just write the resulting matrices.

$$
P_a(2)=\left(
\begin{array}{cc|c}
    0    & 0    & 1 \\ \hline
	ye_2 & 0    & 0 \\ \hline
	ye_1 & ye_2 & 0 \\
	0    & ye_1 & ye_2 \\
\end{array}
\right)
$$ 

Here the first two stripes belongs to $R_1$ and $R_2$ respectively and the last stripe belongs to $R_{12}'$. So we have generators:
 \begin{equation}
 \begin{split}   
   & v_1, v_2, v_3 \in R_3,\\
   &u_1 \in R_1,\\
   & u_2 \in R_2,\\
   &u_1^{12}, \b{u}_1^{12}, u_2^{12}, \b{u}_2^{12} \in R_{12}'.
 \end{split} \tag{*}
 \end{equation} 
  Note that $ye_1 u_i^{12} = \b{u}_i^{12}$ and $ye_2 u_i^{12} = y u_i^{12} - \b{u}_i^{12}$ for $u_i^{12} \in R_{12}'$, $i=1,2$. Then we have the following relations for these generators:
\begin{align*}
&	z' v_1 = ye_2 u^2 + ye_1 u_1^{12}\\
&	z' v_2 = ye_2 u_1^{12} + ye_1 u_2^{12} \\
&	z' v_3 = u^1 + ye_2 u_2^{12} \\
 \intertext{It implies that} 
&	\b{u}_1^{12} = z' v_1 - y u^2\\
&	\b{u}_2^{12} = z' v_1 - y u^2 + z' v_2 - y u_1{12} \\
&	u^1          = z' v_1 - y u^2 + z' v_2 - y u_1{12} + z' v_3 - y u_2^{12} 
\end{align*}
Now we can exclude generators $\b{u}_1^{12}, \b{u}_2^{12}, u^1$. It is important to note that $\b{u}_i^{12}$, $i=1,2$ are annihilated by $x$. And since $u^1 \in R_1$ is also annihilated by $x$, $u^2 \in R_2$ is annihilated by $z$ and $u_1^{12}, u_2^{12} \in R_{12}'$ are annihilated by $x z$ we have the following relations for $v_1, v_2, v_3, u^2, u_1^{12}, u_2^{12}$ (*):
\begin{align*}
&    z u^2 = 0 \\
 &   x z u_1^{12} = 0 \\
  &  x z u_2^{12} = 0 \\
&	x z' v_1 - x y u^2 = 0 \\
&	x z' v_2 - x y u_1^{12} = 0 \\
&	x z' v_3 - x y u_2^{12} = 0 \\
\end{align*}
It gives the following matrix factorization with columns corresponding to $u^2, u_1^{12}, u_2^{12}, v_1, v_2, v_3$, in this order:
{
$$
Q_a(2)=\left(
\begin{array}{cccccc}
	z & 0        & 0        & 0            & 0            & 0     \\
	0     & xz & 0        & 0            & 0            & 0     \\
    0     & 0        & xz & 0	           & 0            & 0     \\
	-xy   & 0        & 0        & xz' & 0            & 0      \\
	0     & -xy      & 0        & 0            & xz' & 0       \\
	0     & 0        & -xy      & 0            & 0            & xz' \\
\end{array}
\right)
$$ 
  
For the other three words with $a_i$, namely $r_3 - \ti{a}_i \sim a_i$, $\ti{a}_i \sim a_i - c_3$, $\ti{a}_i \sim a_i$ we obtain the matrix factorizations by excluding some generators and the appropriate srows and columns:
\begin{itemize}
 \item Excluding $u^2$ from the list of generators (*) and deleting the first row and 1st column from the matrix $Q_a(2)$ we get the matrix factorization for $\ti{a}_i \sim a_i - c_3$.
 \item Excluding $v_3$ from the list of generators (*) and deleting the last row and the last column from the matrix $Q_a(2)$ we get the matrix factorization
  for $r_3 - \ti{a}_i \sim a_i$.
 \item Excluding both $u^2, v_3$ from the list of generators (*) and deleting the first and the last rows and the first and the last columns from the matrix $Q_a(2)$ we get the matrix  factorization for $\ti{a}_i \sim a_i$.
\end{itemize}

Now one can easily see how the matrix factorization $Q_a(i)$ for the word $\fW_a(i) = r_3 - \ti{a}_i \sim a_i - c_3$ looks like for $i>2$.

\section{Generators and relations. Other words}
\label{s5}

  For other modules of the first level the corresponding matrix factorizations are calculated in a similar way. We only present the results for $n=2$,
  since otherwise we obtain too cumbersome matrices.
  
For the word $\fW_b(2) = r_4 - \ti{b}_2 \sim b_2 - c_4$ we have the matrix of correspondences with columns corresponding to $u^1, u_1^{12}, u_2^{12}, v_1, v_2, v_3$:  
{
$$  
Q_b(2)=\left(
\begin{array}{cccccc}
	x         & 0         & 0         & 0            & 0            & 0     \\
	0         & xz  & 0         & 0            & 0            & 0     \\
    0         & 0         & xz  & 0	         & 0            & 0     \\
	-xy       & -xy       & 0         & xz' & 0            & 0      \\
	0         & 0         & xy        & 0            & xz' & 0       \\
    -zy & -zy & zy  & zz' & zz' & zz' \\
\end{array}
\right)
$$ }
  
For the word $\fW_c(2) = r_4 - \ti{c}_2 \sim c_2 - r_3$ we have the matrix of correspondences with columns corresponding to $u^2, u_1^{12}, u_2^{12}, v_3, v_2, v_1$: 
{
$$
Q_c(2)=\left(
\begin{array}{cccccc}
	z & 0        & 0        & 0            & 0            & 0     \\
	0     & xz & 0        & 0            & 0            & 0     \\
    0     & 0        & xz & 0	           & 0            & 0     \\
	0     & 0        & -xy      & xz' & 0            & 0      \\
	0     & -xy      & 0        & 0            & xz' & 0       \\
	-xy   & 0        & 0        & 0            & 0            & xz' \\
\end{array}
\right)
$$ }

For the word $\fW_d(2) = c_4 - \ti{d}_2 \sim d_2 - c_3$ we have the matrix of correspondences with columns corresponding to $u^2, u_1^{12}, u_2^{12}, u_3^{12}, v_4, v_3, v_2, v_1$:   
$$
Q_d(2)=\left(
\begin{array}{cccccccc}
	z & 0        & 0        & 0            & 0            & 0            & 0 & 0    \\
	0     & xz & 0        & 0            & 0            & 0            & 0 & 0    \\
    0     & 0        & xz & 0	           & 0            & 0            & 0 & 0    \\
    0     & 0        & 0        & xz     & 0	          & 0            & 0 & 0    \\
	-x    & 0        & 0        & 0            & xz' & 0            & 0 & 0     \\
	0     & 0        & 0        & -xy          & 0            & xz' & 0 & 0     \\
	0     & 0        & 0        & -xy          & 0            & 0            & xz' & 0   \\
	0     & 0        & -xy      & 0            & 0            & 0            & 0            & xz' \\
\end{array}
\right)
$$ 

For the word $\fW_c'(2) = r_4 - \ti{c}'_2 \sim c'_2 - r_3$ we have the matrix of correspondences with columns corresponding to $u^1, u^2, u_1^{12}, u_2^{12}, v_4, v_3, v_2, v_1$: 
$$
Q_c'(2)=\left(
\begin{array}{cccccccc}
	x   & 0     & 0        & 0        & 0                   & 0            & 0            & 0  		   \\
	0   & z & 0        & 0        & 0                   & 0            & 0            & 0  		   \\
	0   & 0     & xz & 0        & 0                   & 0            & 0            & 0  		   \\
    0   & 0     & 0        & xz & 0	                & 0            & 0            & 0   		  \\
    0   & 0     & 0        & 0        & xzz' & 0            & 0            & 0   		  \\
	0   & 0     & 0        & -xy      & 0                   & xz' & 0            & 0   		  \\
	0   & 0     & -xy      & 0        & 0                   & 0            & xz' & 0    		 \\
	-xy & -xy   & 0        & 0        & -xyz'      & 0            & 0            & xz' \\
\end{array}
\right)
$$ 

For the word $\fW_d'(2) = c_4 - \ti{d}'_2 \sim d'_2 - c_3$ we have the matrix of correspondences with columns corresponding to $u^2, u_1^{12}, u_2^{12}, u_3^{12}, v_5, v_4, v_3, v_2, v_1$, namely $Q_d'(2)$ equals:   
$$
Q'_d(2)
\left(
\begin{array}{ccccccccc}
	z & 0        & 0        & 0            & 0   & 0         & 0            & 0 & 0    \\
	0     & xz & 0        & 0            & 0   & 0         & 0            & 0 & 0    \\
    0     & 0        & xz & 0	           & 0   & 0         & 0            & 0 & 0    \\
    0     & 0        & 0        & xz     & 0	 & 0         & 0            & 0 & 0    \\
	-x    & 0        & 0        & 0            & xz' & 0            & 0 & 0     \\
	-x    & 0        & 0        & 0            & xzz' & -xzz' & 0            & 0 & 0     \\
	0     & 0        & 0        & -xy          & 0   & 0         & xz' & 0 & 0     \\
	0     & 0        & 0        & -xy          & 0  & 0          & 0            & xz' & 0   \\
	0     & 0        & -xy      & 0            & 0  & 0          & 0            & 0            & xz' \\
\end{array}
\right)
$$ 

For the truncated words (without the first or the last letter) we apply the procedure analogous to that described at the end of the preceding section.

In this way we obtain all matrix factorizations of the polynom $F$ corresponding to the modules of the first level.


 \end{document}